\documentstyle[12pt]{article}

 \input{amssym}

\begin{document}

\vspace{3cm }

\centerline{ \large \bf On the cohomology algebra of free loop spaces.}
\vspace{3mm}

\centerline{ \small by}
\vspace{3mm}

\centerline{ {\bf {\sc Bitjong Ndombol and Jean-Claude Thomas}}}

\vspace{3mm}
 \noindent University of  {\sc Dschang}  - Cameroun and  University of {\sc
Angers - CNRS-6093}  France

 D\'epartement de
Math\'ematiques -  2, Bd Lavoisier
 49045 Angers. France.

\vspace{ 0,3 cm}

\rule{16 cm}{0,2mm}
\vspace{ 0,2 cm}

       \noindent {\bf Abstract.} {\small
       Let $X$ be a simply connected space and $\Bbb K$ be any field. The
normalized singular cochains $N^*(X; {\Bbb K})$ admit a natural strongly
homotopy commutative algebra structure, which induces a natural product on
the
Hochschild
       homology $HH_* N^*X$ of the space $X$. We  prove that, endowed
with this product, $HH_*N^*X$ is isomorphic
       to the cohomology algebra of the free
       loop space of $X$ with coefficients in $\Bbb K$. We also show how to
construct  a simpler Hochschild complex which
       allows direct
       computation.}

\vspace{ 0,2 cm}
\rule{16 cm}{0,2mm}

       \vspace{6mm}

       \noindent {\bf Introduction.}

 The classical definition of the normalized Hochschild complex ${\frak
C}_*A$ of an algebra $A$ extends
naturally to a differential graded cochain algebra $(A, d_A)$ over a
field ${\Bbb K}$ of
characteristic
$p \geq 0$ (see  I-$\S 2$). Suppose  $A$ is  augmented   and let  $\rho:
(A\otimes BA, D) \rightarrow (BA, {\bar D})$  denote the
canonical
projection on  the reduced bar construction $BA$. The homology of the chain
complex ${\frak
C}_*A = (A\otimes BA,
D)$  is the Hochschild homology (with coefficients
in $A$) of
$(A,d_A)$ and is
  denoted by  $HH_*A$. If the graded algebra $A$ is commutative it is
well-known that
       $H^\ast B A$ and $HH_*A $ are commutative graded algebras.

  The purpose of the first part of this paper is to construct a product on
${\frak C}_*A$ when
$A$ is not commutative. For this we first embed the category {\bf DA} of
augmented differential graded
algebras into the category {\bf DASH} of strongly homotopy graded algebras
introduced by H.J. Munkholm
in 1974,
\cite{[M]}.
The objects of {\bf DA} and {\bf DASH} are the same but {\bf DASH}$(A,A')$
= {\bf DC}$(BA,BA')$ where
{\bf DC} denotes the category of coaugmented differential graded coalgebras.
Then we consider the
subcategory of $shc$-algebras (strongly homotopy commutative graded
algebras)  of which the objects are
differential graded algebras
$(A,d_A)$ with multiplication
 $m$ belonging to ${\bf DASH} (A\otimes A, A)$ and
 satisfying some natural axioms (see I-$\S3$). In particular, such algebras
have  graded commutative
homology algebras $H(A,d_A)$ and it is well-known that    there is a
natural Hopf algebra
       structure  on $BA$
       such that $H^*BA$ is a commutative Hopf algebra. First we establish :

       \vskip 2mm

       {\bf Theorem 1} {\it Let $(A,d_A)$ be a $shc$-algebra. The
$shc$-structure  induces a natural graded commutative
 algebra structure
on   the Hochschild homology  $HH_*A $. Moreover,
the canonical projection  $\rho$ induces a natural algebra map $\rho_*:
HH_*A\rightarrow H_*BA$.}

\vskip 2mm

        The product structure of $HH_*A$ is given
       explicitly on the normalized Hochschild complex ${\frak C}_*A $.
       Nonetheless, complete
       computations
       are not tractable directly from ${\frak C}_*A$. To overcome this
       difficulty we introduce ($\S 5$) the notion of
       $shc$-equivalence between two $shc$-algebras. If $(A,d_A)$ and
       $(A',d_{A'}) $ are $shc$-equivalent, then the algebra
       $HH_*A$ is isomorphic to
       $ HH_*A'$. A particular case of interest is that of a differential
graded
       algebra $(A,d_A)$ which is   $shc$-equivalent to
a commutative
       differential graded algebra. This is the case, for instance,   if
$(A,d_A)$ is the
       algebra of normalized singular cochains of a space $X$
       with coefficients in $\Bbb K $ whenever

       a) $X$ is a finite dimensional smooth manifold and ${\Bbb K}={\Bbb R}$,
       (example I-$\S3.4$),

       b) $X$ is a connected topological space and ${\Bbb K}={\Bbb Q}$ (example
       I-$\S5.4$),

       c)  $X$ is an $r$-connected finite complex and  ${\Bbb K}$ is of
characteristic $p$,  with $p>
       \frac{\dim X}{r}$ ($r>1$) (example I$\!$I-$\S 4$-3).

One more particular case is that of a $shc$-algebra
       $(A,d_A, \mu_A)$  which is
       $shc$-equivalent to the commutative graded algebra
       $H(A,d_A)$ equipped with   zero differential. This special case is
far from being trivial, as
illustrated by  examples   1,2,3,4, in section I$\!$I-$\S 4$, or the
computations made by K. Kuribayashi
       \cite{[Ku]}.

       When $H^0A={\Bbb K}$, $H^1A=0$ and $dimH^iA<\infty$ for every $i$,  the
$shc$-equivalence class of an
$shc$-algebra
$(A,d_A)$ can be represented by a {\it shc-model}.
       This is a free non-commutative model of the differential graded algebra
$(A,d_A)$ which
is a quotient of $\Omega B(A,d_A)$,  enriched by a structural map. In this
case, a simpler Hochschild
complex is constructed, (see I-$\S6$). This is the main tool for
     the  computations in examples 1, 4 and 5, in section
I$\!$I-$\S 4$.

 \vskip 3mm

       Let $X$ be a topological space and denote simply by $N^*X$ the
${\Bbb K}$-algebra $N^*(X;
{\Bbb K})$ of normalized singular cochains of $X$. By theorem 1, the
natural $shc$-structure on $N^*X$
allows one to define a natural  commutative graded algebra structure on $
HH_*N^*X$. Denote
by $X^{S^1}$  the free loop space of the topological space $X$,  that is
the space of all continuous
maps from the circle into $X$. In 1987, J.D.S. Jones \cite{[J]}
constructed an isomorphism of graded
vector spaces
       $HH_*N^*X \cong H^{*}(X^{S^1}; {\Bbb K} )$.

       Using theorem 1 and the acyclic model theorem for cochain functors,
we prove that Jones'
isomorphism is compatible with the product  on $HH_*N^*
X $ defined in theorem 1 and the usual cup
product on $H^{*}(X^{S^1}; {\Bbb K})$. More precisely:

       \vskip 2mm

       {\bf Theorem 2} {\it Let $X$ be a simply connected space.  There
exists a  natural equivalence
of cochain  complexes ${\frak C}_*(N^*X) \rightarrow C^*(X^{S^1})$ which
induces a
       natural algebra
       isomorphism

\centerline{$HH_*N^*X \cong H^{*}(X^{S^1}; {\Bbb K}).$}

Furthermore the algebra map $\rho_*: HH_*N^*X \rightarrow H^{*}BN^*X$
identifies with
       $ j^*: H^{*}(X^{S^1}; {\Bbb K} )\rightarrow H^{*}(\Omega X; {\Bbb
K})$,  where $j:
\Omega
       X\rightarrow X^{S^1}$
       denotes the inclusion associated to some base point in $X$.}

 This  extends results obtained by M. Vigu\'e-Poirrier,
       \cite{[V]}, N. Dupont-K. Hess \cite{[DH1]} and by S. Halperin and M.
       Vigu\'e, \cite{[HV]}.  The algebraic techniques
       introduced in I-$\S 6$ allow one  to make  computations in the general
       case. For instance in $\S 4$ -examples 5 and 6 of part I$\!$I  we
       examine the cases
       $X= \Sigma {\Bbb C}P^2$,  ${\Bbb K}={\Bbb F}_2$ and $X=G_2$, the
       exceptional Lie
       group, with ${\Bbb K}={\Bbb F}_5$.

       \vskip 2mm

       Theorem 2 is somehow related to the formalism of path-integrals in
       supersymmetric quantum mechanics as explained below.

       Let $X$ be a finite dimensional smooth manifold and let $LX$ be the
space
       of smooth maps from $S^1$ to $X$. This space may be
       given the structure of an infinite dimensional manifold modelled on a
       Fr\'echet space. If $X$ is paracompact, the natural inclusion
       $LX \subset X^{S^1}$ is a homotopy equivalence. On the other hand, if
       $A_{DR}(-)$ denotes the functor ``differential forms'' then
       the ``iterated integral map'' defined by K.T. Chen, \cite{[Ch]},
       $$
       {\frak C}_* A_{DR}(X) \to A_{DR}(LX)$$ is a homomorphism of differential
       graded algebras, which is a quasi-isomorphism when $X$ is
       simply connected. Since $N^*X$ is $shc$-equivalent to $A_{DR}(X)$ there
       exists an
       isomorphism of graded algebras
       $$
       H_{DR} (LX) \cong H^*(LX; {\Bbb R})\,.
       $$
       This is  a de Rham theorem for the infinite dimensional manifold $LX$.

       In the light of the result of Gezler, Jones and Petrack,
\cite{[G-J-P]}, one
       may ask : {\it Is  the natural equivalence ${\frak C}_*(N^*X)
\rightarrow C^*(X^{S^1})$
       a $shm$-map or  a homomorphism of $A_\infty
       $-algebras in the sense of Stasheff, \cite{[Sta]}? }

       \vskip 3mm

       We would like to thank N. Dupont, E.H Idrissi and M. Vigu\'e
       for their interest and for valuable suggestions during
       the preparation of this paper. We are also very grateful to J.
Stasheff, G. Powell and to the referee
for helpful comments.

\vspace{6mm}

       \centerline{\large{ \bf Part I: Algebraic setting.}}
       \vspace{3mm}

       Throughout the paper, the ground field $\Bbb K$ of characteristic
$p\geq 0$ is
       fixed. We use the  Kronecker convention: an object with lower negative
       graduation has upper non-negative graduation.

       \vspace{3mm}

       \noindent {\bf 1. Review of the bar and cobar constructions}

       \vskip 3mm
       \noindent {\bf 1.1} Recall that $ \bf DA$ (resp. $\bf DC$ resp. $\bf
DM$)
       denotes the
       category of augmented,
       differential graded algebras (resp. coaugmented differential graded
       coalgebras, resp.
       differential graded modules).
       An object  $A \in\mbox{Obj \bf DA} $ is a graded $\Bbb K$-vector
space  $A = \{A^k\}_{k\geq 0}$, equipped with structure
       $$
       d_A: A^k \to A^{k+1} \,,  m: \oplus_{k+l=n} A^l\otimes A^k\to A^n\,,
       \eta_A: {\Bbb K}\to A \,, \epsilon_A: A\to {\Bbb K}
$$
and an exact sequence
$0\rightarrow IA \stackrel{i_A }{\rightarrow}A\stackrel{\epsilon_A
       }{\rightarrow}{\Bbb K}
       \longrightarrow 0$.  Similarly,  an object $C\in \mbox{Obj\bf DC} $
is a graded $\Bbb K$-vector space  $C =
\{C^k\}_{k\geq
0}$, with structure
$$
 d_C: C^k \to C^{k+1}\,,  \triangle: C^n\to \oplus_{k+l=n} C^l\otimes C^k\,,
       \epsilon_C: C\to {\Bbb K} \,,  \eta_C: {\Bbb K}\to C
$$
and an exact sequence $ 0\rightarrow {\Bbb K} \stackrel{\eta_C
       }{\rightarrow}C\stackrel{j_C
       }{\rightarrow}JC
       \longrightarrow 0 $.

       The morphisms in $\bf DA$ (resp. $\bf DC$) respect  the whole
structure and
       are called homomorphisms of $DG$-algebras (resp.
       homomorphisms of $DG$-coalgebras).

       \vspace{3mm}

       \noindent {\bf 1.2} We denote by
       $
       {\underline B}: \mbox {\bf DA }\rightarrow \mbox {\bf DC} \,, \mbox{
       resp. (} {\underline \Omega}:
       \mbox {\bf DC} \rightarrow \mbox {\bf DA} \mbox{)}
       $
       the bar construction, (resp. the cobar construction) and by
       $
       B: \mbox {\bf DA} \rightarrow \mbox {\bf DC} \mbox{ resp. (} { \Omega}:
       \mbox {\bf DC} \rightarrow \mbox {\bf DA} \mbox{)}
       $
       the reduced bar construction (resp. the reduced cobar construction)
       \cite{[A]}, \cite{[ML]}, \cite{[L]}, \cite{[M]},
       \cite{[FHT]}. The functors  $B$ and   $\Omega$ (resp. ${\underline
B}$ and
${\underline
       \Omega}$) are adjoint functors to each other. We denote by
       $ \omega: \mbox {\bf DA} (\Omega C, A) \longrightarrow \mbox {\bf DC}(C,
       BA)$ the natural bijection.
       As ${\Bbb K}$-graded vector spaces
       $$
       {\underline B}A = T(A), \quad \quad BA = T(IA),\quad {\underline
\Omega C} =
       T(C)\quad \mbox{and} \quad \Omega C = T(JC)
       $$
       and we denote  by $
       [a_1\vert a_2\vert ...\vert a_k] \mbox{ (resp. } \langle c_1\vert
c_2\vert ....\vert c_l\rangle \mbox{)} $ the standard  generators of $B^kA$
or $\underline{B}^k A$ (resp.
$\Omega
       ^l C$ or $\underline{\Omega }^l C$ ) of degree
       $\sum_{i=1}^k \deg a_i-k$ (resp.
       $\sum_{i=1}^l \deg a_i+l$).

       The natural inclusion ${B}A \rightarrow \underline{B} A$ (resp. the
       quotient ${\underline \Omega}C \rightarrow \Omega C$) induces  an
       isomorphism in (co)homology.

       \vspace{3mm}
\noindent {\bf 1.3} The homomorphism $\alpha_A\in \mbox{\bf DA}(\Omega BA, A)$
        corresponding, by adjunction, to the identity $id_{BA}$ is a
natural   homotopy
       equivalence \cite{[M]},
       \cite{[FHT]}. More precisely, $\iota_A\in \mbox{\bf DM}(A, \Omega
BA)$, defined by
       $$
       \iota_A(a)=
       \left\{
       \begin{array}{rl}
       \eta_A\epsilon_A(a) &\mbox{ if } \epsilon_A(a)\neq 0\\
       \langle [a]\rangle &\mbox{ if} \quad \epsilon_A(a)= 0
       \end{array}
       \right.
       $$ satisfies $\alpha_A\circ \iota_A = id_A$ and $id_{\Omega
BA}-\iota_A \circ
       \alpha_A =
       d_{\Omega BA}\circ h +h\circ d_{\Omega BA}$, for some chain homotopy
 $h:\Omega BA \to \Omega BA$ such
       that $\alpha_A\circ h = 0$, $h\circ \iota_A = 0$, $h^2 = 0$.

       \vspace{3mm}
       \noindent {\bf 1.4} Recall that $f, g \in \mbox{\bf DA}(A, A')$ are
       homotopic in $\mbox{\bf DA}$
       if there exists a
       linear map $h: A\to A'$ such that $f - g = d_{A'}\circ h + h\circ
d_A$ and
       $h(xy) =
       h(x)g(y) + (-1)^{\vert x\vert}f(x)h(y)$ with $x$, $y\in A$.
       Observe that $\iota_A \circ \alpha_A \simeq id_{\Omega BA}$ in
$\mbox{\bf
       DM}$ but {\em not}
       in $\mbox{\bf DA}$.

\vspace{3mm}
       \noindent {\bf 2. Hochschild homology of a differential graded cochain
       algebra.}

       \vskip 3mm
       \noindent {\bf 2.1} Let $(A, d_A)$ be a cochain algebra:
       $$
       A = \{A^i\}_{i\geq 0}\,, \quad d_A: A^k \rightarrow A^{k+1} \,,
\qquad x\in
       A^i\,, \quad
       \vert x\vert := i\,.
       $$
 Let $d_{BA}$ denote the differential of the bar
       construction ${\underline B}A$ or $BA$. The tensor product
       $(A,d_A) \otimes ({\underline B}A, d_{BA})$ (resp. $(A,d_A)\otimes
(BA, d_{BA}$) is then a
differential
       module whose differential is denoted
       $d_{A\otimes BA}$ in
       both cases. The Hochschild
       differential, denoted by $D$, is defined by:
       $$
       Da_0[a_1\vert
       ....\vert a_n] = (d_0 - d_n)a_0[a_1\vert
       ....\vert a_n] + d_{A \otimes BA} a_0[a_1\vert
       ....\vert a_n]
       $$
       where
       $d_0 a_0[a_1\vert
       ....\vert a_n] = (-1)^{\vert a_0\vert} +
       a_0a_1[a_2\vert ...\vert ..\vert a_n]$ and
       $d_na_0[a_1\vert
       ....\vert a_n] =(-1)^{(\vert a_n\vert +1) (\vert a_0\vert
       +...+\vert a_{n-1}\vert +n-1)} \\
       a_na_0\otimes[a_2\vert ...\vert ..\vert a_{n-1}]$.
By definition,
       $${\frak C} _* A := (A\otimes {B}A, D) \,,
       \mbox{( resp. }
       {\underline{\frak C}} _* A := (A\otimes {\underline B}A, D)\mbox{)}
       \,.$$
is the
       {\it (normalized)
       Hochschild complex} (resp. {\it (un-normalized) Hochschild complex}) of
       $(A, d_A)$ and
       $$
       HH_*A=H{\frak C}_*A = H{\underline{\frak C}}_*A
       $$
       is the {\it Hochschild homology } of the cochain algebra $(A, d_A)$. One
       should notice here that $\frak{C}_*A $ is concentrated
       in non-negative total upper degrees.
       In particular  $ HH_*A $ is concentrated in non-negative upper degrees.

       \vspace{3mm}
       \noindent {\bf 2.2} Let $\Sigma _{n,m}$ be the set of
$(n,m)$-shuffles. Consider the shuffle
map,
\cite{[L]} 4.2.1,
       $sh: \underline{\frak C}_* A\otimes \underline{\frak C}_* A
       \longrightarrow \underline{\frak C _*
       }(A\otimes
       A)$ defined by:
       $$
sh(a_0[ a_1 \vert a_2\vert ... \vert a_n] \otimes  b_0[ b_1
       \vert b_2\vert ... \vert
       b_m] ) = (-1)^t \sum_{\sigma \in \Sigma_{n,m}}
       (-1)^{\epsilon (\sigma)}a_0\otimes b_0[
       c_{\sigma(1)}\vert ...\vert c_{\sigma (m+n)}]
       $$
       where
       $t = \vert b_0\vert (\vert a_0\vert +...+\vert a_{n}\vert)$,
       $
       \quad \quad
       c_{\sigma (i)} = \left \{ \begin{array}{ll}
       a_{\sigma (i)}\otimes 1, &1\leq i \leq n \\
       1\otimes b_{\sigma (i-n)}, &n+1\leq i \leq n+m
       \end{array}
       \right.
       $\\
       and $\epsilon (\sigma) = \sum (\vert c_{\sigma (i)}\vert-1)( \vert
c_{\sigma
       (m+j)}\vert-1)$, summed over all pairs
       $(i, m+j)$ such that $\sigma (m+j) <\sigma (i)$.
Clearly, $sh$ induces a chain map  $sh : {\frak C}_*
       A\otimes {\frak C}_* A \longrightarrow {\frak C
       _* }(A\otimes
       A)$ and a map of differential
       coalgebras
       $ sh : BA\otimes BA\rightarrow B(A\otimes A)$.

       \vspace{3mm}
       \noindent {\bf 3 $shc$-algebras.}

       \vskip 3mm
       \noindent {\bf 3.1} A {\it strongly homotopy commutative algebra}
({\it $shc$-algebra}
       for short)
       is a triple $(A, d_A, \mu_A)$ with
       $(A, d_A)\in \mbox{Obj\bf DA}$ and $\mu_A \in \mbox{\bf DA}(\Omega
       B(A\otimes
       A), \Omega BA)$ satisfying
       \begin{enumerate}
       \item $\alpha_A\circ \mu_A \circ \iota_{A\otimes A} = m_A$,  where
$m_A$ is
the product in
$A$;
       \item $\alpha_A\circ \mu_A \circ \Omega B(id_A\otimes \eta_A)\circ
\iota_A =
\alpha_A\circ \mu_A\circ
       \Omega B(id_A\otimes\eta_A)\circ \iota_A
       =id_A$;
       \item $\mu\circ \Omega B(\alpha_A\otimes id_A)\circ \Omega
B(\mu\otimes id_A)\circ
       \chi_{(A\otimes A)\otimes A} \simeq
       \mu\circ \Omega B(id_A\otimes \alpha_A)\circ
       \Omega B(id_A\otimes \mu)\circ \chi_{A\otimes( A\otimes A)}$ in
$\mbox{\bf
       DA}$;
       \item $\mu\circ \Omega BT \simeq \mu$ in $\mbox{\bf DA}$,
       \end{enumerate}
where  $T$ denotes  the interchange map $ T(x\otimes y ) = (-1)^{|x|\, |y|} y
\otimes x$. The two  natural homomorphisms of DG-algebras, defined in
\cite{[M]}-2.2,
$$
\Omega B(\Omega B(A\otimes A)\otimes A)\stackrel{\chi_{(A\otimes A)\otimes
A }}{\longleftarrow}
\Omega B(A\otimes A\otimes A)
        \stackrel{\chi_{A\otimes (A\otimes A)}} {\longrightarrow }
       \Omega B(A\otimes \Omega B(A\otimes A))
 $$
 satisfy: $
       \alpha_{(A\otimes A)\otimes A}\circ \chi_{(A\otimes A)\otimes A}
       =\alpha_{A\otimes A\otimes A} =
       \alpha_{A\otimes (A\otimes A)}\circ \chi_{A\otimes( A\otimes A)}\,.$

       In particular, if $(A, d_A, \mu_A)$ is a $shc$-algebra,  then
conditions $1$
       and $4$
       in the definition imply that $H^*(A)$ is a commutative graded algebra.

       \vspace{3mm}
       \noindent {\bf 3.2} Consider $A$ and $A'$ in $\mbox{Obj\bf DA}$. The
linear map $f\in
       \mbox{\bf DM}(A, A')$
       is  said to be a
       {\it $shc$-map } from $(A, d_A , \mu_A)$ to $(A',d_A', \mu_A')$ if
there exists
       ${\underline f}\in \mbox{\bf DA}(\Omega BA, \Omega BA')$ such that:
       \begin{enumerate}
       \item $\alpha_A\circ {\underline f}\circ \iota_A = f$;
       \item $\alpha_{A'}\circ {\underline f}\circ \eta_{\Omega BA} =
\eta_{A'}$;
       \item ${\underline f}\circ \mu_A \simeq \mu_A'\circ \underline
{f\otimes f}$ in {\bf DA}.
       \end{enumerate}
where $\underline {f\otimes f}$, is defined in  \cite{[M]}-2.2. Moreover,
if $f\in \mbox{\bf DA}(A,
A')$,  then
$f$ is a {\it strict $shc$-map}. Observe that  a
       $shc$-map is a $shm$-map in the sense of \cite{[M]}. Obviously, if
$A$ and $A'$ are commutative
graded algebras, then any $f\in
       \mbox{\bf DA}(A,
       A')$ is a strict $shc$-map.

       \vspace{3mm}
       \noindent {\bf 3.3 Proposition.} {\it Suppose that  $(A, d_A,
\mu_A)$ is a
       $shc$-algebra,  then $BA$ is a Hopf algebra which is associative up to
homotopy, and $H BA$ is a
commutative
       Hopf algebra.
       }

       \vskip 3mm
       \noindent {\bf Proof.} Take  $\omega $ as defined in 1.2. Observe that:

\centerline{ $\alpha_A\circ \mu_A\in
       \mbox{\bf DA}(\Omega B(A\otimes A), A) \mbox{ and } sh\in \mbox{\bf
DC}(BA\otimes BA, B(A\otimes
A))$}
\noindent  hence $\nu_A=\omega (\alpha_A\circ
       \mu_A)\in
       \mbox{\bf DC}(B(A\otimes A), BA)$ and the composite
       $\nu_A\circ sh \in \mbox{\bf DC} (BA\otimes BA, BA)$ defines a Hopf
algebra
structure       on $BA$. Moreover,
$\nu_A\circ T
\simeq
\nu_A$. Indeed $\nu_A\circ T =
\omega (\alpha_A\circ
       \mu_A)\circ BT = \omega (\alpha_A\circ \mu_A\circ T) \simeq \omega
(\alpha_A\circ \mu_A) = \nu_A$, since $\omega$ preserves homotopy.

       \vspace{3mm}
       \noindent {\bf 3.4 Example.} Let $X$ be a finite-dimensional smooth
       manifold and let $\cal U$ be an open cover of $X$. We denote
       by $A_{DR}(X)$ the de Rham complex of $X$ and by $N^*(X; {\Bbb R})$ the
       cochain complex of normalized singular cochains with
       coefficients in $ \Bbb R$.  The natural homomorphism of cochain
complexes
       $
       A_{DR}(X) \to N^*(X; {\Bbb R})$, $ \quad \omega \mapsto
\int_{\,}\omega \,,
       $
       is a $shm$-map by \cite{[B-G]}. Consider the \v Cech-de Rham bicomplex,
       $
       N^{**}({\cal U}; A_{DR}(X)) = \bigoplus _{p,q \geq 0} N^p({\cal U};
       A^q_{DR}(X)) $ with  the usual differential and with product defined by
       $
       (\alpha\cup
       \beta )_{\vert _{U_{i_0}\cap U_{i_1}\cap ...\cap U_{i_{p+p'}}}}
       = (-1)^{p'} \alpha _{\vert _{U_{i_0}\cap U_{i_1}\cap ...\cap
       U_{i_{p}}}} \cup \beta _{\vert _{U_{i_p}\cap U_{i_{p+1}}\cap ...\cap
       U_{i_{p+p'}}}}.
       $
Now from \cite{[M]} proposition 4.7 we see that $N^*(X; {\Bbb R})$ is a
       $shc$-algebra. In the same way, so is the
       \v Cech-de Rham complex,
       $
       N^* ({\cal U}; A_{DR}(X))$. Moreover the two natural maps
       $$
       \begin{array}{rl}
       N^*(X; {\Bbb R}) \to N^* ({\cal U}; A_{DR}(X)),  &
       \, c \in N^p(X; {\Bbb R})\mapsto (c_{\vert _{U_i}}) \in N^p ({\cal U};
       A^0_{DR}(X))\\
       A_{DR}(X) \to N^* ({\cal U}; A_{DR}(X)), &
      \,  \omega \in A^q(X) \mapsto (\omega_{\vert _{U_i}}) \in N^0 ({\cal U};
       A^q_{DR}(X))
       \end{array}
       $$
       are strict $shc$-maps and, for a good choice of the open cover $\cal
U$,
       these are quasi-isomorphisms.

       \vspace{3mm}
       \noindent
       {\bf 4 Proof of theorem 1}

       \vspace{3mm}
       \noindent {\bf 4.1}
       Let $(A, d_A, \epsilon _A)$ be an augmented cochain algebra. Since
the cochain map
       $\iota _A: A\rightarrow \Omega BA$, (1.3)
       is not a homomorphism of $DG$-algebras, it does not induce a
homomorphism
       of Hochschild complexes.
       Nonetheless ${\frak C _* }(\alpha_A):{\frak C _* }(\Omega BA)
       \rightarrow {\frak C _* }(A)$
       is a surjective quasi-isomorphism. Therefore there exists a chain map
       $$
s_A:{\frak C _* }A \rightarrow {\frak C _* }\Omega BA
$$
such that ${\frak C _* }\alpha_A\circ s_A =Id_{{\frak C}_* A} $ and
       $s_A\circ  {\frak C}_* \alpha_A$ is  chain homotopic to
       $id_{{\frak C} _* (\Omega BA)}$.

Now suppose that  $(A, d_A, \mu_A)$ is a $shc$-algebra. We define
       a homomorphism of
       $DG$-modules
       $$
       \Phi: \underline{\frak C} _* A\otimes
       \underline{\frak C }_* A \rightarrow \underline {\frak C} _* A\,,
       \quad \quad \Phi = \underline{\frak C} _* \alpha_A\circ
       \underline{\frak C} _* \mu\circ
       s_{A \otimes A} \circ sh.
       $$
where $s_{A\otimes A}$ denotes  the linear section of $\alpha_{A\otimes A}$,
as defined above.
 Notice that $\Phi$ induces a homomorphism of
       cochain complexes
       ${\frak C}_* A \otimes {\frak C}_* A \rightarrow {\frak C} _* A$, also
       denoted by $\Phi $. Since $Hs_{A\otimes A} =
       (H\alpha_{A\otimes A})^{-1}$, the homomorphisms
       $H\Phi$ does not depend on the
       choice of the section $s_{A \otimes A}$.
       Actually, precomposing $H_*(\bar \Phi)$ by the $K\ddot{u}nneth$
isomorphism
       $HH_*A\otimes
       HH_*A \rightarrow H_*({\frak C _* }A\otimes
       {\frak C _* }A)$ yields a multiplication
       $$
       \Phi_*: HH_*A\otimes HH_*A \rightarrow HH_*A \,.
       $$
       The associativity property of $\Phi _*$ is a direct consequence of the
       associativity of the
       shuffle map and of $H\mu_A$ together with the fact that
       the morphisms $\iota_{- }$ and $\chi_{- }$ induce the identity in
homology. It
       results from axiom 4 in the definition of a $shc$-algebra that $H\eta_A:
       k\rightarrow HH_*A$ is a unit. Axiom 2 in the definition of
        a  shc-structure and
       naturality imply the commutativity of $HH_* A$.

       \vspace{3mm}

       \noindent {\bf 4.2 Lemma.} { \it If $\varphi: A \rightarrow B$ is a
strict
       $shc$-map
       then $HH_*\varphi $ is a homomorphism of graded algebras. Moreover, if
       $\varphi$ is  a quasi-isomorphism of DG-algebras then
       $HH_*\varphi $
       is an isomorphism.}
       \vskip 2mm

       \noindent {\bf Proof.} If $\varphi$ is a homomorphism of DG-algebras
then
       $HH_*\varphi$ is a well-defined linear map, which is an
       isomorphism when $\varphi$ is a quasi-isomorphism. The fact that
$HH_*\varphi$
       preserves multiplications follows directly from the definition of a
$shc$-map
       together with the naturality
       of the constructions.

       \vspace{3mm}

       \noindent {\bf 4.3} Finally,  observe that the diagram
$$
       \begin{array}{cccccccccc}
       {\frak C}_* A \otimes {\frak C}_* A &
       &\stackrel{sh}{\rightarrow} & {\frak C}_* (A \otimes A)
       &\stackrel{{\frak C}_* (\alpha_{A\otimes A})}{\leftarrow} & {\frak C}_*
       \Omega B(A\otimes A)
       &\stackrel{{\frak C}_* \mu_A}{\rightarrow} & {\frak C}_* \Omega BA
       & \stackrel{{\frak C}_* \alpha_A}{\rightarrow} & {\frak C}_* A \\

       \Big\downarrow {\rho \otimes \rho}&&
       &\Big\downarrow \rho &
       &\Big\downarrow \rho &
       &\Big\downarrow \rho &
       & \Big\downarrow \rho \\

       BA\otimes BA &
       & \stackrel{sh}{\rightarrow}& B(A\otimes A)
       &\stackrel{ B\alpha_{A\otimes A}}{\leftarrow} & B\Omega B(A\otimes A)
       & \stackrel{B\mu_A}{\rightarrow}& B\Omega BA
       &\stackrel{B\alpha_A}{\rightarrow}& BA\\
       \end{array}
$$
       commutes. Now, using adjunction, we have that $\nu_A\circ
       B\alpha_{A\otimes A} =
       B\alpha_A\circ B\mu_A$. This ends the proof of theorem 1.

       \vspace{0,4 cm}

\vfill{\eject}

       \noindent {\bf 5 $shc$-equivalence and $shc$-formality.}

       \vspace{3 mm}
       \noindent {\bf 5.1 } Two $shc$-algebras $A$ and $A'$ are
{\it $shc$-equivalent}
       if there exists a sequence of strict $shc$-maps,
       $$
       A \leftarrow A_1 ... \rightarrow A_2 \leftarrow ... \rightarrow A'
       $$
       which are quasi-isomorphisms. This implies,  in particular, that
$HH_* A$ and
$HH_* A'$
are isomorphic as graded algebras.

 \vspace{4 mm}

       \noindent {\bf 5.2 } The $shc$-algebra $A$ is said to be {\it
$shc$-commutative}
if $(A, d_A , \mu _A)$ is
       $shc$-equivalent to a commutative differential graded
       algebra $(A', d_{A'})$.

       The $shc$-algebra $A$ is said to be {\it $shc$-formal} if $A$ is
$shc$-equivalent to the
       commutative differential graded
       algebra $H(A,d_A)$ equipped with zero differential; in this case
$HH_* A\cong HH_* (H(A))$ can be
computed using   differential forms,  as
established by  Hochschild, Kostant and Rosenberg, \cite{ [H-K-R]}.

 \vspace{3mm}
       \noindent {\bf 5.3} As proven in \cite{[H]}, any commutative
       graded algebra $ A$ admits  a free commutative model ${\cal M}_A$
which is  of the form
$(\Lambda X, d)\rightarrow
       A $,
       where
       $\Lambda X = E(X^{odd})
       \otimes
       P(X^{even})$ and $E$ means exterior algebra and $P$ polynomial
algebra. Therefore, if $(A,d_A)$
is $shc$-commutative, then $HH_*A= HH_* {\cal M}_A$ and   one deduces easily
from \cite{[HV]} and \cite{[ML]}-$\S 8$-2.3 that
there exists a commutative graded differential algebra of the form   $ (\Lambda
       X\otimes \Gamma sX, d')$,  where $\Gamma$
       stands for the free algebra of divided powers,   such that
       $$H^{\ast }(\Lambda
       X\otimes \Gamma sX, d') \cong HH_* A \mbox{ and } H^\ast (\Gamma sX,
       \bar d') \cong H^\ast BA
       $$
       as  graded algebras.

       \vspace{3mm}
       \noindent {\bf 5.4 Example.} Let ${\Bbb K}$ be a field of
characteristic
       zero and let $X$ be a connected space.   A simple generalization
       of example 3.4 with the polynomial de Rham forms,
       $A_{PL}(X)$, playing the role of $A_{DR}(X)$ and a simplicial \v Cech-de
       Rham bicomplex $N^{**}( X ; A_{PL}(X))$
       (see \cite{[B-G]} or \cite{[FHT2]}-$\S 10$ for more details) proves
that
the $shc$-algebra
       $N^*(X; {\Bbb K})$  is $shc$-commutative. If there exist
quasi-isomorphisms of $DG$-algebras $ A_{PL}(X) \rightarrow (A,d_A)
\leftarrow H^*(X, {\Bbb K})$
       where   $A$ is
       a commutative graded algebra, then  the space $X$ is  said to be ${\Bbb
K}$-formal and  $N^*(X;{\Bbb
       K})$ is $shc$-formal. This is no longer true for a field
       of positive characteristic, as shown by  example I$\!$I-4.4 or
I$\!$I-4.5.
However, it follows  from example 3.4 and \cite{[DGMS]} that,
       if
       $X$ is a compact K\"alher manifold or a Riemannian symmetric space,
and  if ${\Bbb K}= {\Bbb R}$,
then
       $HH_* N^*X  = HH^* H_{DR}(X)$ as an algebra.

       \vspace{0,4 cm}
       \noindent {\bf 6 A Small Hochschild complex for a $shc$-algebra.}

\vspace{3mm}
       \noindent {\bf 6.1} Let $(A, d_A, \mu_A)$ be a $shc$-algebra. The
composite $BA\otimes BA\stackrel {sh} \rightarrow B(A\otimes A) \stackrel
{\nu_A} \rightarrow BA$, where $\Omega \nu_A = \mu_A$, defines the product
(denoted $\star$) considered in 3.3. In particular, since $\Omega (BA\otimes
BA) = (T((s^{-1}B^+A\otimes {\Bbb K})\oplus ({\Bbb K}\otimes s^{-1}B^+A)\oplus
s^{-1}(B^+A\otimes B^+A)), D)$, the composite
$$
\mu_0: (s^{-1}B^+A\otimes {\Bbb K})\oplus ({\Bbb K}\otimes s^{-1}B^+A)\oplus
s^{-1}(B^+A\otimes B^+A) \stackrel {s^{-1}sh} \rightarrow s^{-1}B^+(A\otimes A)
\stackrel {s^{-1}\nu_A} \rightarrow s^{-1}B^+A
$$
called the ''linear part'' of
$\mu$ satisfies

$\mu_0(s^{-1}[a_1\vert ...\vert a_i]\otimes 1) = s^{-1}[a_1\vert ...\vert
a_i]\,, \quad \mu_0(1\otimes
s^{-1}[a_1\vert ...\vert a_i]) = s^{-1}[a_1\vert ...\vert a_i]$

$\mu_0(s^{-1}[a_1\vert ...\vert a_i]\otimes  [b_1\vert ...\vert b_j]) =
s^{-1}([a_1\vert ...\vert a_i]\star [b_1\vert ...\vert b_j]).$

       \vspace{3mm}

 \noindent {\bf 6.2} Let $(A,d_A)$ be an augmented differential graded
algebra and suppose we have an
algebra quasi-isomorphism $(TU, D)\rightarrow (A, d_A)$. One example of
this situation is
 $\alpha_A: \Omega BA\rightarrow (A,
d_A)$) where $U = s^{-1}B^+A$,
 $D=D_1+D_2$
       with $D_1U\subset U$ while
       $D_2U\subset T^{\geq 2}U$. Here, $D_1$
is, up to a shift of degrees,
$d_{BA}$  the differential of the bar
       construction,  and $D_2 u = \sum_i\langle a_i\vert
       a'_i\rangle$ where
       $\bar \triangle u = \sum_i a_i\otimes a'_i$ is the reduced coproduct in
       $BA$.

 \vspace{3mm}
       \noindent {\bf 6.3}  Let $(A,d_A)$ be an augmented graded algebra
and suppose that:

1) $H^0(A, d_A) = {\Bbb K} $ and that $ H^1(A, d_A) = 0$

2) There is a differential $D$ on $TU$ such that $D= D_1+D_2$ with $D_1 U
\subset U $ and $D_2 U \subset T ^{\geq 2} U$.

3) we have chosen a quasi-isomorphism $(TU,D) \to (A,d_A)$ of DG-algebras,

   Set $U = kerD_1 \oplus S$, $kerD_1 = D_1S \oplus V$. Clearly, $V \cong
H^\ast (U, D_1)$.
The  decreasing filtration $F^pTU = \oplus_{i\geq p} T^iU$ yields a first
quadrant spectral sequence which
converges to $H^\ast TU
\cong  H^\ast A.$ Let $I$ denote the ideal generated by $S$ and  $DS$ in
$TU$;   the ideal $I$ is
acyclic and hence the
projection $p: (TU, D) \rightarrow (TU/I, {\bar D})$ is a quasi-isomorphism.

   From the decomposition $U = V \oplus S \oplus D_1S$, we deduce an
isomorphism of graded algebras
$ TU/I \cong TV $ which defines a differential $d_V$ on $(TV, d_V)$ and a
surjective
quasi-isomorphism
$$
 p_V: (TU, D) \rightarrow (TV, d_V) \,.
$$

 Moreover, there exists a homomorphism of $DG$-algebras $ \varphi _V : (TV,
d_V) \to (TU, D)$ such that
$p_V \circ \varphi _V= id$. Therefore $\varphi _V$ is a quasi-isomorphism.
Observe that when $(TU, D)
=  \Omega BA$,
$$ \Omega B A \stackrel{p_V} \rightarrow (TV, d_V)
       \stackrel{\varphi _V}\rightarrow \Omega B A \,,$$
and

       a) $V= \{V^i\}_{i \geq 2}$ and it is isomorphic to $H(U, D_1) \cong
s^{-1} H^+(BA)$,

       b) the linear part of $d_V$ is zero while its quadratic part is, up
to a shift of degree, the
  comultiplication of $H^\ast BA$.

The differential graded algebra $(TV, d_V)$ is called a free minimal model
for $(A, d_A)$. Obviously,
this model depends on the direct factor $V$ in $U$. It is easy to see that
another choice produces a
minimal model $\varphi _{V'} : (TV', d_{V'})  \rightarrow (TU, D)$ and an
isomomorphism $ \varphi  :
(TV, d_V) \to (TV', d_{V'})$ such that $
       \varphi_{V'} \circ \varphi \simeq \varphi _V$. For this
       reason, one can speak of the minimal model of a space. (See also
\cite{[FHT]}).

 \vspace{3mm}
       \noindent {\bf 6.4}  Let $(A,d_A, \mu_A)$ be an augmented
$shc$-algebra and
assume that
$H^0(A, d_A) = {\Bbb K} $ and that $
       H^1(A, d_A) = 0$. Write $\Omega (BA\otimes BA) = (T{\hat U}, d_ {\hat
U})$. We denote by $\mu _V : (T\hat V , \hat d) \to (TV,d)$ the composite
$$ (T \hat V, d_{\hat
V} )
   \stackrel{ \varphi _{\hat V}}{\rightarrow }  \Omega (BA\otimes BA)
   \stackrel{\mu_A\circ sh}{\longrightarrow}  \Omega BA
 \stackrel{ p_V}{\rightarrow }  (TV, d_V)
$$
The triple
$(TV, d_V, \mu _V)$ is called a {\it $shc$-minimal model for $(A, d_A, \mu
_A)$}.

 \vspace{3mm}
       \noindent {\bf 6.5} In the remainder
 of this  section we  establish  some of
the main properties of the  $shc$-minimal model of the $shc$-algebra
$(A, d_A, \mu_A)$. First observe that,  by  6.2, 6.3 and 6.4, we have a
direct sum decomposition:
 $$
\begin{array}{rl}
\hat V &\cong s^{-1} (H^+BA \otimes {\Bbb K}) \oplus s^{-1} ({\Bbb K} \otimes
H^+ BA) \oplus s^{-1} (H^+BA \otimes H^+BA)
\\ &\cong  s^{-1}(sV \otimes {\Bbb
K}) \oplus s^{-1}({\Bbb K} \otimes sV) \oplus s^{-1} (sV \otimes sV) \,.
\end{array}
$$
Therefore  we write:
$\quad
\hat V := V  \oplus  W \oplus V \# W \,,$ with $V := s^{-1}(sV \otimes {\Bbb
K})$, $W :=  s^{-1}({\Bbb K} \otimes sV) $, $  V \# W := s^{-1} (sV \otimes
sV) $ and
      $ v\#w := s^{-1}(sv\otimes sw) $.

Define $\tilde \Psi:\Omega (BA\otimes
BA) = (T((s^{-1}B^+A\otimes {\Bbb K})\oplus ({\Bbb K}\otimes s^{-1}B^+A)\oplus
s^{-1}(B^+A\otimes B^+A)), D) \rightarrow \Omega BA\otimes \Omega BA$ by:
$\tilde \Psi(s^{-1}[a_1\vert ...\vert a_i]\otimes 1) = s^{-1}[a_1\vert ...\vert
a_i]\otimes 1$, $\tilde \Psi(1\otimes s^{-1}[a_1\vert ...\vert a_i]) =
1\otimes s^{-1}[a_1\vert
...\vert a_i]$, $\tilde \Psi(s^{-1}[a_1\vert ...\vert a_i]\otimes
[b_1\vert ...\vert b_j] )=
0$. One can  check  easily that   $\tilde \Psi$ commutes with the
differentials. Moreover, the commutativity of the diagram
$$
\begin{array}{cccc}
& \Omega B(A\otimes A)&
 \stackrel{\alpha_{A\otimes A}}{\longrightarrow} &A\otimes A \\
 & \Omega sh \Big\uparrow &
&\Big\uparrow \alpha_{A}\otimes \alpha_{A} \\
&
 \Omega (BA\otimes BA) &\stackrel{\tilde \Psi}{\longrightarrow} & \Omega
BA\otimes
\Omega BA\\
\end{array}
$$
\noindent implies that $\tilde \Psi$ is a
surjective quasi-isomorphism. Set
$$
\Psi_V=(p_V\otimes p_V)\circ \tilde \Psi\circ  \varphi _{\hat
V}: (T{\hat V}, d_{\hat V})\rightarrow (TV, d_V)\otimes (TV, d_V)\,,
$$
hence $\Psi_V$ is also  a surjective quasi-isomorphism and since $p_{\hat
V}\circ \varphi _{\hat
V} = id$,  $\Psi_V$ satisfies:

$\Psi_Vv = v\otimes 1$, \quad $\Psi_Vw = 1\otimes w$,  \quad $\Psi_Vv\#w = 0$,
\quad $v\in V$ and $w\in W$.

\vspace{3mm}
       \noindent {\bf 6.6 Proposition.} {\it With the previous identification
then

1) $\mu_V $ identifies $V$ and $W$ with  $V$,

2) if $\mu_0 : \hat V= V\oplus W \oplus V \# W \to V $ denotes the linear
part of $\mu_V$, then
$\mu_0(v\# w) = 0 \mbox{ if and only if }  sv\star sw =0 $ where $\star $
is the product in $H^+BA
\simeq sV$,

3) the differential $d_{\hat V}$ is completely determined by the fact that
$\Psi_V: (T{\hat V},
d_{\hat V})\rightarrow (TV, d_V)\otimes (TV, d_V)$ is a surjective
quasi-isomorphism,

4) if  either $d_V v=0$ or $d_V w=0$,  then $d_{\hat V} (v \# w )$ is
determined
by the Hirsch-type formulae described below.}

\vspace{3mm}
\noindent{\bf Proof} Parts 1 and 2 of the proposition are  a direct
consequence  of the definitions.  The rest of this section is devoted to
the proof  of parts 3 and 4.

 The degree $-1$ identification $ \quad
       V\# W \leftrightarrow s^{-1}(sV\otimes sW)\,, \quad
       v\#w \leftrightarrow s^{-1}(sv\otimes sw) $
defines a  linear map of degree $-1$:
$V\otimes W\rightarrow T(V\oplus W\oplus V\# W)$ which we wish  to extend to
$T(V)\otimes T(W)$. Set
$$
\begin{array}{crl}
\mbox{a) }&  (v_1...v_k)\# w &= \sum^{k}_{i=1}(-1)^{a_i}v_1..v_{i-1}(v_i\#
w)v_{i+1}...v_k \\
                            & & v_i\in V \,, w\in W\,, a_i = \vert v_1\vert
+...+\vert v_{i-1}\vert
+\vert w\vert (\vert v_{i+1}\vert
       +....+\vert v_k\vert) \\
\mbox{b) }&  v\#( w_1...w_l) &= \sum^{l}_{j=1}(-1)^{b_j}w_1..w_{j-1}(v\#
       w_j)w_{j+1}...w_l \\
    & & v\in V \,, w_j\in W \,,  b_j = (\vert v\vert
+1)(\vert
       w_1\vert +....+\vert w_{j-1}\vert ) \\
\mbox{c) }&  (v_1...v_k)\# w_1...w_l &=
       \sum^{k}_{j=1}(-1)^{c_j}w_1..w_{j-1}(v_1..v_k\# w_j)w_{j+1}...w_l\\
&&v\in V \quad w_j\in W \,, c_j =
       (\vert (v_1...v_k)\vert +1)(\vert w_1\vert +....+\vert w_{j-1}\vert )

\end{array}
$$

       Formula c) means that we first expand $\#$ with respect to the first
argument  and
       then with respect to the second argument. If one  reverses this order,
i.e. we
expand
       first with respect to the second argument and
       then  with respect to the first argument,  one
       obtains the formula
$$
\begin{array}{crl}
\mbox{d) }& (v_1...v_k)\#' w_1...w_l &:=
       \sum^{k}_{j=1}(-1)^{d_i}v_1..v_{i-1}(v_i\# w_1...w_l)v_{i+1}...v_k\\

       && v_i\in V \quad w_j\in W $, $ d_i =
       (\vert (v_1\vert +...\vert v_{i-1} +(\vert v_{i+1}\vert ..\vert
       v_k\vert)(\vert w_1\vert +....+\vert w_l\vert )
\end{array}
$$
Observe that, in general, $(v_1...v_k)\# w_1...w_l \neq (v_1...v_k)\#'
w_1...w_l$ and that formula a)
can be read as $ (v_1...v_k)\#' w$.  To avoid confusion we keep the
notation used in formula a).

We define a derivation of degree $1$, $D_0: T(V\oplus W\oplus V\# W)
\rightarrow T(V\oplus W\oplus
       V\# W)$ by:
$$
       D_0v = d_Vv\,, \quad D_0w =
       d_Ww\,,\quad
       D_0(v\#w) = v \cdot w - (-1) ^{\vert v\vert \, \vert w\vert} w\cdot v -
       D_0v\#w -(-1)^{\vert
       v\vert}v\# D_0w
$$
       where $ v\in V$, $ w \in W$.

\vskip 3mm
\noindent {\bf 6.7 Lemma} {\it For $v \in V \,,\,  w\in W \,,\,  a \in
TV\,, \, b\in TW $,
       \begin{enumerate}
       \item $D_0(a\#'w) = a\cdot w - (-1) ^{\vert a\vert \, \vert w\vert}
       w\cdot a - (D_0)a\#'w -(-1)^{\vert
       a\vert}a\#' D_0w$,

       \item $D_0(v\#b) = v  b - (-1) ^{\vert v\vert \, \vert b\vert}
       b v -
       D_0v\#b -(-1)^{\vert
       v\vert}v\# D_0b$.
       \end{enumerate}
       }
 \vskip 3mm
       \noindent {\bf Proof.} We proceed by the induction on word-length.
First observe that
equalities 1 and 2 hold when $a\in V$ and $b\in W$. Suppose the first
equality true when $a\in T^pV$
and consider
$a=a_0a_1\in
       T^{p+1}V$ with
       $a_0\in T^pV$
       and $a_1\in V$. Then by  the induction hypothesis,
$$
       \begin{array}{rl}
       D_0(a\#'w) = & (-1)^{\vert a_0\vert} D_0a_0.(a_1\#'w) +
       a_0(a_1w-(-1)^{\vert a_1\vert \vert w\vert}
       wa_1) -a_0(D_0a_1\#'w-(-1)^{\vert a_1\vert} a_1\#'D_0w)\\
       &+ (-1)^{\vert a_1\vert \vert w\vert} a_0wa_1 -
       (-1)^{\vert a_0\vert \vert w\vert +\vert a_1\vert \vert w\vert }
wa_0a_1\\
       & - (-1)^{\vert a_1\vert \vert w\vert}(D_0a_0\#'w)a_1 -
       (-1)^{\vert a_1\vert \vert w\vert +\vert a_0\vert}
       (a_0\#'D_0w)a_1 \\
&  + (-1)^{\vert a_1\vert +1)\vert w\vert +\vert a_0\vert +\vert w\vert
+1}(a_0\#'w).D_0a_1
       \end{array}
$$
 so that $ D_0(a\#w) = aw - (-1)^{\vert a\vert \vert w\vert}wa -
D_0a\#w -
       (-1)^{\vert a\vert}a\#D_0w.$

       The second equality is obtained in the same way.

\vskip 3mm
       \noindent {\bf 6.8} {\it End of the proof of proposition 6.5.}
Recall that $\Psi_V$ is a
surjective quasi-isomorphism. The end of the proof of part 3 follows
directly from the uniqueness,
up to isomorphism, of the minimal
       model (6.3). Moreover,  we deduce from lemma 6.7  that
$$
       D_0^2(v\#w) = D_0v\#' D_0w - D_0v\# D_0w,
$$
 which is not zero in general. Thus   $d_{\hat V}$ does not
coincide with
$D_0$. Nonetheless we have
$d_{\hat V} (v \# w)  = D_0(v \# w)$, if $v$ or $w$ is a cocycle in $V$.

 \vspace{3mm}
       \noindent {\bf 6.9 Proposition.} {\it Let ${\cal M}_A :=(TV, d_V,
\mu_V)$ be a
$shc$-model of the $shc$-algebra $(A, d_A,
\mu_A)$. There exists a product on the Hochschild homology $HH_*{\cal M}_A $
which coincides, up to an isomorphism, with the
product defined on $HH_*A$.}

\vskip 3mm

\noindent {\bf Proof} Consider the following diagram

$$
 \begin{array}{ccccccc}
 {\frak C}_*A\otimes {\frak C}_*A &
 \stackrel{sh}{\longrightarrow} & {\frak C}_*(A\otimes A) &
\stackrel{{\frak C}_*\alpha_{A\otimes A}}{\longleftarrow}
&{\frak C}_*\Omega B(A\otimes A)& \stackrel{{\frak C}_*\mu_A\circ \alpha_A}
{\longrightarrow} & {\frak C}_*A \\
\Big\uparrow &    &\Big\uparrow &   &\Big\uparrow & & \Big\uparrow \\
{\frak C}_*TV\otimes {\frak C}_*TV &\stackrel{sh}{\longrightarrow} & {\frak
C}_*(TV\otimes TV) &
\stackrel{{\frak C}_*\Psi}{\longleftarrow}
&{\frak C}_*T{\hat V})& \stackrel{{\frak C}_*\mu_V}{\longrightarrow} &
{\frak C}_*TV
\end{array}
$$
\noindent where the vertical arrows  are respectively from left,
 ${\frak C}_*(\alpha_A \circ \varphi _V ) \otimes  {\frak C}_*(\alpha_A
\circ \varphi _V )$,
  ${\frak C}_*(\alpha_A \otimes \alpha _A \circ \varphi_V\otimes \varphi_V)$,
 ${\frak C}_*( \varphi_{\hat V})$
and ${\frak C}_*(\alpha_A\circ \varphi_V)$. The  two left
hand squares commute, while the right hand square
 commutes up to homotopy.

If $s_{\hat V}$ denotes a linear section of ${\frak C}_*\Psi$, one defines
the product $\Phi_V$ on ${\frak C}_*TV$ by
$\Phi_V = {\frak C}_*\mu_V\circ s_{\hat V}\circ sh$.

\vspace{1 cm}

\centerline {\large{ \bf Part II: From  algebra to topology.}}

\vspace{3mm}
\noindent {\bf 1. Hochschild homology of a space.}

\vspace{3 mm}

\noindent {\bf 1.1} Let $X$ be a topological space. We denote by  $C^*X$
(resp. $ N^*X$)
the algebra of un-normalized singular cochains (rep. of normalized
singular cochains) on
$X$. Since the  inclusion $N^*X\rightarrow C^*X$ is an algebra map and
induces an
isomorphism in
(co)homology, we  define  the Hochschild homology of $X$ as the graded
vector space
$$
  HH_*N^*X \cong HH_*C^*X\,.
$$

 From \cite{[M]}1.2, we know that  there exists a  commutative \underline
{natural} diagram
in
${\bf DM}$
$$
\begin{array}{cccc}
& \Omega B(N^*X\otimes N^*X)&
 \stackrel{\underline {AW}}{\longrightarrow} &\Omega BN^*(X\times X) \\
 & \alpha_{N^*X\otimes N^*X} \Big\downarrow &
&\Big\downarrow \alpha_{N^\ast (X\times X)} \\
&
 N^*X\otimes N^*X&\stackrel{AW}{\longrightarrow} & N^*(X\times X)\\
\end{array}
$$ where  $AW$ denotes the normalized Alexander-Whitney map.

\vspace{3mm}
\noindent{ \bf 1.2  Proposition.} {\it   Let $\triangle$ be the topological
diagonal map. If $X$ is path connected,  the
natural $shc$-structural map of
$N^\ast X$  defined by
$
\mu_X = \Omega BN^*\triangle \circ \underline {AW}
$
  induces  a  natural graded commutative algebra structure on
 $HH_*X$.}

\vspace{3mm}
\noindent {\bf 1.3 Example.}  A space $X$ is called  {\it ${\Bbb
K}$-$shc$-formal} if $N^\ast
X$ is $shc$-formal. Spheres and complex projective spaces are $shc$-formal
for any field.  If $X$ is a
simply connected space
 which is
${\Bbb K}$-$shc$-formal, then the multiplicative structure on $H^*(X^{S^1};
{\Bbb K})$  is  completely
determined by the graded algebra $H^*(X; {\Bbb K})$. See \cite{[Id]} and
examples 4.1, 4.2 below.

\vspace{3mm}
\noindent {\bf 2. The acyclic model theorem for cochain  functors}

\vspace{3 mm}
 The proof of theorem 2  relies heavily on the acyclic model theorem for
cochains functors. For the convenience of the reader we
recall some definitions here.

\vspace{4 mm}

\noindent{\bf 2.1}  Let us denote  the category of
non-negatively graded vector spaces by ${\bf DM^*}$. Let ${\cal A}$
be a category with models ${\cal M}$.  Recall that a functor $F: {\cal A}
\rightarrow{\bf DM^*}$

a) {\it admits a unit} if,  for each object $A$ in ${\cal A}$,  there exists a
linear map $ \eta_A: {\Bbb K} \rightarrow FA$ such that $d\circ \eta_A = 0$;

b)  {\it is acyclic on the models} if,  for any object $M$ in ${\cal M}$,
there
exists a linear map $\epsilon_M: FM\rightarrow {\Bbb K} $
such that $\epsilon_{M}\circ \eta_M \simeq Id_{\Bbb K}$ and $ \eta_M\circ
\epsilon_M
\simeq Id_{FM}$.

c) is {\it corepresentable on  the models} if there
exists a natural transformation
$\kappa: {\hat F}\rightarrow F$ such that $\kappa \circ \xi = id_F$, where
$:\hat F$ denotes the  contravariant functor
$$
{\hat F}: {\cal A}\rightarrow {\bf DM^*}
\,, \quad {\hat F}(A) =
\prod_{M\in {\cal M}} F(M)\times {\cal A}(M, A)\,.
$$
and $ \xi $ the  natural transformation:
$ \xi_A: F(A) \rightarrow {\hat F}(A) $, $ a\mapsto \{F(f)(a)\}_{M
\in {\cal M}\,,f\in
{\cal
A}(M, A)}\,, a \in A\,. $

\vspace{3 mm}
For instance, the functor $X\mapsto  C^*X$ is corepresentable on
the standard simplexes $\Delta ^n$ in ${\bf Top}$.

\vspace{3mm}
\noindent {\bf 2.2 Theorem} (\cite{[B-G]})
{\it Let ${\cal A}$ be a category with models ${\cal M}$ and
$F_1, F_2: {\cal A}\rightarrow {\bf DM}^*$
two contravariant functors with units. If $F_1$ is acyclic and $F_2$ is
corepresentable on the models, then :\\
1) there exists a natural transformation $\tau: F_1 \rightarrow F_2$ which
preserves the units,\\
2) any   two such natural transformations are naturally homotopic.}

\vspace{3mm}
\noindent {\bf 3. Proof of theorem 2.}

\vspace{3 mm}
\noindent {\bf 3.1} Consider the simplicial set $K$ defined as follows:
$K(n) = {\Bbb Z}/(n+1){\Bbb Z}, $ and, if $\overline {k}^{n}$ denotes an
element in ${\Bbb Z}/n{\Bbb Z}$,
the face maps  $d_i: K(n)\rightarrow K(n-1)$ with $0\leq i\leq n-1$ and the
degeneracy maps  $s_j:  K(n)\rightarrow K(n+1)$
with
$0\leq i\leq n$
are:
 $$
d_i \overline { k}^{n+1}=\left\{
\begin{array}{rl}
\overline {k}^{n} & \mbox{ if } k\leq i\\
\overline {k-1}^{n} & \mbox{ if } k>i
\end{array}
\right.
\qquad
 s_j\overline {k}^{n+1}=
       \left \{
       \begin{array}{rl}
        \overline { k}^{n+2}  &\mbox{ if } k\leq i\\
       \overline {k+1}^{n+2} &\mbox{ if }  k>i .
       \end{array}
       \right.
$$
 and $ d_n \overline {k}^{n+1}  = \overline {k}^{n}$.  Consider also the
simplicial set $P$ defined by $P(n) = \overline
{0}^{n+1}\in {\Bbb
Z}/(n+1){\Bbb Z}$ with obvious face and degeneracy maps. If $\Sigma$ is a
simplicial set,  we denote as usual by $\vert
\Sigma \vert $ its geometric realization, \cite{[ML]}.
By \cite{[BF]} (proposition 1.4),  $\vert K\vert$ is
homeomorphic to the circle $S^1$ while  $\vert P \vert$
is a point.

\vspace{3mm}
\noindent {\bf 3.2} To any  topological space  $X$ one
can associate the cosimplicial topological spaces $\underline X$,
$\underline Y$ and
$\underline P$ defined by:
$
\underline{X}(n)  = Map (K(n), X) =
\underbrace {X\times....\times X}_{(n+1) \mbox{-times}}$,   ${\underline
P}(n) = Map(P(n), X) = X\quad n\geq 0 $,
and  if $* \in X$,
$   \underline{Y}(n)  = \{*\}\times
\underbrace {X\times....\times X}_{(n) \mbox{-times}}\,, \quad n\geq 0 $.

Observe  that $Y(n) =  \{f\in Map (K(n), X) \mbox{ such that }
f(\overline{0}^{n+1})=*
\}.$

The coface and codegeneracy maps for the  cosimplicial  spaces $\underline
X$ and the sub-cosimplicial space
$\underline{Y} \subset \underline X$ are:
$$  \begin{array}{ll}
d_i(x_0,x_1,...,x_n) &= (x_0,x_1,...,x_i,x_i, ... x_n) \,, 0\leq i \leq n\\
d_{n+1}(x_0,x_1,...,x_n) &= (x_0,x_1,...,x_n, x_0)\\
s_j(x_0,x_1,...,x_n)&= (x_0,x_1,...x_j,x_{j+2},...,x_n)\,, 0\leq j \leq n.
\end{array}
$$
We have a sequence of obvious cosimplicial maps
$\underline Y\to \underline X\to \underline P$.

\vspace{3mm}
\noindent {\bf 3.3} Write  ${\bf Top}$ (resp. ${\bf Costop}$) for  the
category of topological spaces (resp. of cosimplicial topological spaces).
There is a covariant functor $\vert \vert . \vert \vert: {\bf
Costop}\rightarrow
{\bf Top}$ called the geometric realization
$$
 \vert \vert \underline {Z}\vert \vert =
{\bf Costop}(\underline{\triangle}, \underline{Z})
\subset \prod_{n\geq 0} {\bf Top}(\triangle^n,{\underline Z}(n))\,,
$$
where $
 \vert \vert \underline {Z}\vert \vert$ is
equipped with the
topology induced by this inclusion. Here $\underline{\triangle}$ denotes
the cosimplicial space defined by
$\underline{\triangle}(n)
= \triangle^n$
 with the usual coface and codegeneracy maps $\delta_i$ and $\sigma_j$
respectively.

If $\Sigma$ is a simplicial set and $T$ a topological space then $\underline Z
= Map(\Sigma, T)$ is a cosimplicial topological space. We recall the following
duality result due to Bott and Segal,     \cite{[BS]} (proposition 5.1).
\vskip 3mm
\noindent {\bf 3.4 Proposition} {\it There is a  homeomorphism:
$$ \vert \vert T^{\Sigma}\vert \vert = \vert \vert \underline Z \vert \vert =
 {\bf Costop}(\underline{\triangle}, \underline{Z}) \cong  {\bf Top}(\vert
\Sigma
\vert , T) = T^{\vert
\Sigma
\vert} \,.$$
}
 With this notation introduced in  3.2 this   yields that
$$ \vert \vert \underline X \vert \vert \cong  {\bf Top}(\vert
\Sigma
\vert , X) = X^{S^1} \mbox{ and } \vert \vert \underline P \vert \vert \cong
X$$ where $X^{S^1}$ is the free loop space of $X$.

\vspace{3mm}
\noindent \noindent {\bf 3.5 Proposition} {\it  $\vert
\vert \underline Y \vert \vert \cong \Omega X$.}

\vskip 3mm
\noindent {\bf Proof} From 3.4 one has
$\vert \vert \underline X \vert \vert \stackrel{G}{\rightarrow }
{\bf Top}(\vert
K
\vert , X) \stackrel{F}{\rightarrow }\vert \vert \underline X \vert \vert$
where $F$ and $G$ are inverse
homeomorphisms.  We denote  by
$\star_0$  the $0$-simplex of
$\triangle_0$ and by
$[\overline {0}^1, \star_0]$ the base point of $\vert K\vert \simeq S^1$.
If $f\in \vert \vert \underline Y\vert \vert = {\bf
Costop}(\underline{\triangle}, \underline{Y})$,
then $f=\{f_n\}_{n\geq 0}$ with
$f_n:\triangle^n\rightarrow
\underline Y(n)$ compatible with the coface and codegeneracy maps  and
satisfying $f_n(s_0^n\overline {0}^1)=\star$, $ n\geq 0$ and where
$s_0^0=id$. By the definition of $G$, $G(f)([\overline {0}^1,
\star_0])=f(\star_0)(\overline {
0}^1)=\star.$ That is,  $G(f) \in {\bf Top}(\vert K
\vert , X)$ preserves the base point. Let $f=\{f_n\}_{n\geq 0}\in {\bf
Top}(\vert
K
\vert , X)$ preserving the base point, that is $f([\overline {0}^1,
\star_0])=\star$. We have $F(f) = \{g_n\}_{n\geq 0}$,  with
$g_n:\triangle^n\rightarrow
\underline X(n)$ compatible with the coface and codegeneracy maps.
By the definition of $F$, $g_n(t)(x) = f_n([t,x])$, $t\in \triangle^n$ and
$x\in
K(n)$. For any  $t\in \triangle^n$, $(t, s_0^n\overline {0}^1)$ is
equivalent to
$(\sigma_0^nt, \overline {0}^1)=(\star_0, \overline {0}^1)$ where the
$\sigma_j$ are the
codegeneracy maps of $\triangle^n$. Thus $g_n(t, s_0^n\overline {0}^1) =
F(f)([t,
s_0^n\overline {0}^1])=\star$ and $F(f) \in {\bf Costop}(\underline{\triangle},
\underline{Y})$.

\vspace{3mm}
\noindent  {\bf 3.6} The naturality of the above constructions yields
the commutative diagram
$$
\begin{array}{rlcccc}
 & \vert \vert \underline Y\vert \vert& \stackrel{ }{\longrightarrow} & \vert
\vert \underline X\vert \vert& \stackrel{ }{\longrightarrow}& \vert \vert
\underline P\vert \vert\\
 & \quad \Big\downarrow   &  & \Big\downarrow  &  &\Big\downarrow
 \\
 &\Omega X &\stackrel{ }{\longrightarrow}& X^{S^1}&\longrightarrow &X \\
\end{array}
$$ in which the vertical arrows are homeomorphisms. In particular the top
line is a Serre fibration.

 If $\underline Z$ is any  cosimplicial topological space,
 then $C^*\underline {Z}$ is a simplicial cochain complex. We define  $Tot
C^*\underline {Z}$ by $(Tot C^*\underline {Z})_n =
\oplus_{p-q=n} C^q\underline{Z}(p)$
with differential $Dx = \sum^p_{i=1}(-1)^iC^\ast(d_i)+(-1)^p\delta x$, where
$x\in
C^\ast \underline Z(p)$, the
 $d_i$ are the coface operators  and $\delta$ is the internal differential of
$C^\ast \underline Z(p)$.
Observe that $ \underline {Z } \mapsto Tot C^*{\underline Z}$ and
$\underline {Z } \mapsto C^*\vert \vert {\underline
Z}\vert
\vert$ are contravariant functors from ${\bf Costop}$ to ${\bf DA}$.

There exists a natural transformation, \cite{[BS]} (corollary 5.3) and
 \cite{[J]} (proof of theorem 4.1 and lemma 6.3):
$$
\psi_{\underline Z}: Tot C^*\underline {Z} \rightarrow C^*\vert \vert
\underline Z\vert \vert$$ such
that $\psi_{\underline X}$ induces an isomorphism in
(co)homology when $X$ is  simply connected. Since it is obvious that
$\psi_{\underline P}: C^*X \rightarrow C^*X$ is the identity map, we have
the commutative diagram
in $\bf DA$
$$
\begin{array}{rlcc}
 &C^*X&\stackrel{ id}{\longrightarrow} & C^*X\\
& \quad \Big\downarrow   &  & \Big\downarrow\\

& Tot C^*\underline {X}&\stackrel{ \psi_{\underline
X}}{\longrightarrow}&C^*\vert \vert \underline X\vert \vert\\

&\quad \Big\downarrow   &  & \Big\downarrow\\

& Tot C^*\underline {Y}&\stackrel{ \psi_{\underline Y}}{\longrightarrow}&
C^*\vert \vert \underline Y\vert \vert . \\
\end{array}
$$
 As pointed out by J. Jones in \cite{[J]} (6. proofs of theorems A and
B), the iteration of the Alexander-Whitney natural transformation $AW:
C^*X\otimes C^*X \rightarrow C^*(X\times X)$ yields a natural transformation
$$ \theta_X: {\frak C}_*C^*X\rightarrow Tot C^*\underline {X}$$ inducing an
isomorphism in (co)homology.

It is straightforward to check that:

-  $\theta_X$ restricts to the identity  on $C^*X$ ,

- $\theta_X$ induces $\theta'_X: BC^*X\rightarrow Tot C^*\underline {Y}$

-  there is a commutative diagram

$$
\begin{array}{rlcccc}
 & N^*X & \stackrel{ i_X}{\longrightarrow}& C^*X &\stackrel{
id}{\longrightarrow} & C^*X\\

&\Big\downarrow   &  & \Big\downarrow& & \Big\downarrow\\

(**) \qquad  &{\frak C}_*N^*X & \stackrel{ \theta_X\circ
j_X}{\longrightarrow}&Tot C^*\underline
{X}&\stackrel{ \psi_{\underline X}}{\longrightarrow}&C^*\vert \vert
\underline X\vert \vert\\

&\rho \Big\downarrow   &  & \Big\downarrow & &\Big\downarrow\\

&BN^*X& \stackrel{ \theta'_X\circ j'_X}{\longrightarrow}&Tot C^*\underline {Y}&
\stackrel{ \psi_{\underline Y}}{\longrightarrow}&C^*\vert \vert \underline
Y\vert
\vert \\
\end{array}
$$
where $i_X:N^*X\to C^*X$ is the inclusion inducing $j_X:{\frak C}_*N^*X\to
{\frak C}_*C^*X$ and $j'_X:BN^*X\to BC^*X$ respectively. The maps $i_X$ and
$\Theta_X = \psi_{\underline
X}\circ \theta_X\circ j_X$ induce isomorphisms in (co)homology.

Set $\Theta'_X = \psi_{\underline
Y}\circ \theta'_X\circ j'_X$, then

\vspace{3mm}
\noindent  {\bf 3.7 Proposition} {\it The map $\Theta'_X$ induces an
isomorphism in (co)homology.}

\vskip 3mm
\noindent {\bf Proof} The map $\rho: {\frak C}_*N^*X \rightarrow BN^*X$ factors
through the diagram
$$
\begin{array}{rlcc}
 & {\frak C}_*N^*X & \stackrel{ \rho}{\longrightarrow}& BN^*X \\
&\rho_1 \Big\downarrow   &  & \rho_3 \Big\uparrow \\
&{\frak C}_*N^*X\otimes_{N^*X}B(N^*X, N^*X)& \stackrel{
\rho_2}{\longrightarrow}& BN^*X\otimes B(N^*X,
N^*X)\\
\end{array}
$$
where  $B(N^*X, N^*X)$ is the bar construction with coefficients, \cite{[ML]},
$\rho_1$ is the inclusion, $\rho_2$ is the usual  isomorphism, $\rho_3$ is
a quasi-isomorphism since
$B(N^*X, N^*X)$ is  a semi-free resolution of the field ${\Bbb K}$ as
$N^*X$-module,
\cite{[FHT1]} (lemma 4.3).  Therefore,  $H^*({\frak
C}_*N^*X\otimes_{N^*X}B(N^*X, N^*X)) = Tor_{N^*X}(
{\frak C}_*N^*X, {\Bbb K})
\cong  H^*(BN^*X)$.

In view of diagram $(**)$, we  also obtain,  \cite{[EM]},  that
$$
H^*(C^*\vert \vert \underline Y\vert \vert )  \simeq Tor_{C^*X}( C^*\vert
\vert \underline X\vert \vert, {\Bbb K})\,,
$$
and that the linear map $H^*(\Theta'_X): H^*BN^* X \to  H^*(C^*\vert \vert
\underline  Y\vert \vert )  $ coincides with
$Tor _{i_X} (
\Theta _X ; {\Bbb K}) : Tor_{N^*X}( {\frak C}_*N^*X, {\Bbb K})\rightarrow
Tor_{C^*X}( C^*\vert \vert \underline X\vert
\vert, {\Bbb K})$. The latter is an  isomorphism, since $i_X $ and
$\Theta_X$ are quasi-isomorphisms,  \cite{[FHT1]}
(proposition 2.3).

Most  of the rest of this section is devoted to  the proof of:

\vspace{4 mm}
\noindent {\bf 3.8  Proposition} {\it If $X$ is simply connected, the
 natural
chain equivalence of
cochain complexes
$$
 \Theta_X:  {\frak C}_* N^*X \longrightarrow C^*\vert \vert  \underline {
X}\vert \vert
$$
 induces an isomorphism of graded algebras in (co)homology.}

\vskip 2mm
We could not prove
each assertion completely without going into too much detail, so we only
prove the main one i.e.  the
commutativity of the following diagram
$$
\begin{array}{cccc}
& HH_*X\otimes HH_*X &
 \stackrel{H(\Theta_X)\otimes H(\Theta_X)}{\longrightarrow} & HC^* \vert \vert
\underline
{X}\vert \vert  \otimes HC^* \vert \vert  \underline {X}\vert \vert \\
(*) & \Phi_* \Big\downarrow &
&\Big\downarrow \cup \\
&
HH_*X &\stackrel{ H(\Theta_X)}{\longrightarrow} & HC^* \vert \vert  \underline
{X}\vert \vert \\
\end{array}
$$
where $ \Phi_*$ is the product considered  in 1.4-(a) and  $\cup$ is the
usual cup product.

 In order to establish the commutativity of the diagram
$(*)$,   consider the natural Alexander-Whitney transformation, \cite{[ML]},
Chap.V$\!$I$\!$I$\!$I-8.1 and 8.6: $$
AW: {\frak C}_*(N^*X\otimes N^*X)  \rightarrow {\frak C}_*N^*X\otimes
{\frak C}_*N^*X
$$
which satisfies
$ AW\circ sh=id$ and $H(sh\circ AW) = id.$

>From this fact and the definition of $s_{N^*X \otimes N^*X}$, the commutativity
 of the above diagram  is a consequence  of the commutativity up to
homotopy of the
following
diagram:
$$
\begin{array}{lcccrl}

\quad {\frak C}_* \Omega B(N^*X\otimes N^*X) &\stackrel{  {\frak
C}_*\alpha_{N^*X\otimes N^*X} }{\to}&
{\frak C} _\ast (N^*X\otimes N^*X) &  \stackrel{AW}{\to} {\frak C}
_*N^*X\otimes {\frak C}_*N^*X\\

&&& \qquad \qquad  \downarrow    \Theta_X
\otimes
\Theta_X \\

{\frak C} _*(\alpha_{N^\ast X}\circ \mu_X)  \Big\downarrow && &C^* \vert
\vert  \underline {X}\vert \vert
 \otimes C^* \vert \vert \underline {X}\vert\vert \\

 &&
&\Big\downarrow \cup \\

\quad  \qquad \qquad {\frak C}_* N^*X &&\stackrel{
\Theta_X}{\longrightarrow}&  C^* \vert \vert
\underline {X}\vert \vert \\
\end{array}
$$
 To establish the commutativity up to homotopy, of this diagram, we
consider the contravariant functors  $\quad F_1,  F_2:
{\bf Top}  \rightarrow {\bf DM}^* \,,$
$$
F_1(X) = {\frak C}_*\Omega B(N^*X\otimes N^*X) \,,\quad  F_2(X) = C^*
\vert \vert  \underline {X}\vert \vert.
$$
  It results from lemma 3.9 below  and from  theorem 2.2 that the two
 natural transformations
$\beta_X = \cup \circ (\Theta_X \otimes \Theta_X)\circ
AW\circ {\frak C}_*\alpha_{N^*X\otimes N^*X}$
and $ \gamma_X = \Theta_X\circ {\frak C}_*(\alpha_{N^\ast X}\circ \mu_X)$
are chain
homotopic and thus the above diagram is commutative up to homotopy.

Indeed, the functor $F_1$ admits a natural unit (i.e. a linear map $\eta _X
: {\Bbb K} \to F_{1}(X)
$ such that $d \circ \eta =0$).  The composition of the
natural inclusions $N^0\vert\vert  {\underline X}\vert \vert  \rightarrow
C^0\vert \vert {\underline
X}\vert \vert  \rightarrow C^*\vert \vert  {\underline X}\vert \vert  $
defines a unit on $F_2$.
By a
tedious verification, one checks  that the natural transformations $\beta$ and
$\gamma$
respect these units.

\vspace{3mm}

\noindent {\bf 3.9 Lemma.}{ \it The functor $F_1$ (resp. $F_2$)  is acyclic
(resp. corepresentable) on
a set of models in ${\bf Top}$.}
\vskip 3mm
\noindent {\bf Proof.}  Let  ${\large \vee}_0^pZ$ be
the wedge of $(p+1)$ copies of $Z$.
Consider the topological spaces:

$O^n =
{\large \vee}_{p\geq
0}({\large \vee}_0^{p+1}(\triangle ^n\times\triangle ^p))$ topologized
with the weak
topology,

 $T_X^n = \{f=\{f_p\}_{p\geq 0}| f_p: \triangle ^n\times\triangle ^p\to
{\underline X}(p)\mbox{ such that } f_p(\star_p)=(f_0(\star_0), ..,
f_0(\star_0))\}$,

$\triangle ^n$ (resp. $\triangle ^n\times\triangle ^p$) pointed by $\star$
(resp.
$\star_p$).

Since $\triangle^n\times \triangle^p$ is contractible for $n,  p\geq
0$, so is
$O^n$ for $n\geq 0$ and $F_1$ is acyclic on the models $O^n$.
To prove the corepresentability of $F_2$ on these models, we establish the
corepresentability of $F(X) =  C^n\vert {\underline X}\vert$,  for each
$n\geq 0$.
Consider
$$
 {\hat F}(X) = \prod_{m \geq 0} \prod_ {f\in {\bf Top}(O^m,
X)}C^n\vert \vert
{\underline O^m}\vert \vert \times \{f\}\,.
$$
The natural transformation
$\xi: F\rightarrow {\hat F}$ is defined by:
$$
\xi (w) =
(C^n\vert {\underline f}\vert w)_f \,, \quad  w\in
C^n\vert {\underline X}\vert
\,,  \quad f\in {\bf Top}(O^m, X), m\geq 0\,
$$
 where ${\underline f}$ denotes the cosimplicial map associated to $f$.

 Let $\sigma \in C_n\vert \vert  {\underline X}\vert \vert$ be a singular
simplex.
Recall that
$\vert \vert {\underline
X}\vert \vert = {\bf Costop}({\underline \triangle},  {\underline X})\subset
\prod_{p\geq
0}{\bf Top}(\triangle ^p, {\underline X}(p))$.

 Observe that there is an
obvious  bijection
$$\vartheta:  {\bf Top}(\triangle ^n, \vert \vert
{\underline
X}\vert \vert )\rightarrow  {\bf Costop}(\triangle ^n\times {\underline
\triangle} ,
{\underline X})$$  where $\triangle ^n\times {\underline
\triangle}$ is the cosimplicial space defined by $(\triangle ^n\times
{\underline
\triangle})(p) = \triangle ^n\times \triangle^p$. Moreover ${\bf
Costop}(\triangle ^n\times {\underline \triangle} ,
{\underline X})$ is a subspace of $T_X^n$. Indeed,  if $\tau
\in {\bf Costop}( \triangle ^n\times {\underline \triangle}, {\underline
X})$, then $\tau
= \{\tau_p\}_{p\geq 0}$ with each $\tau_p: \triangle ^n\times \triangle^p \to
{\underline X}(p)$ compatible with the first coface map $d_0$.

For any $n\geq 0$, there is a natural
bijection
$$ \tau ^X_n: T^n_X
\rightarrow {\bf Top}(O^n, X)\,.$$ defined by $\tau^X_n(g) = {\large
\vee}_{p\geq
0}({\large \vee}_{i=0}^{p+1}g_{p, i}).$

 Thus $\sigma$ determines
a
unique element ${\tilde \sigma} = \tau ^X_n\circ \vartheta (\sigma)$ in
${\bf Top}(O^n,
X)$.\\
Finally observe that, if $id$ is the identity map of $ O^n$, then
$(\tau^{O^n}_n)^{-1}(id)\in {\bf CosTop}(\triangle^n\times {\underline
\triangle}, {\underline O^n}).$ Set
$\iota_n = (\vartheta)^{-1}\circ (\tau^{O^n}_n)^{-1}(id) \in C_n\vert \vert
{\underline O^n}\vert \vert $.\\
Let $(a_{m, f})\in {\hat F}(X)$ and $\sigma \in C_n\vert \vert {\underline
X}\vert \vert$.
The formula
$(\kappa_X(a_{m, f}))(\sigma) = a_{n, \tilde \sigma}(\iota_n)$ defines a
natural
transformation $\kappa: {\hat F}\rightarrow F$ such that $\kappa \circ \xi
= id$, as can be seen from the diagram:

$$
\begin{array}{cccc}
& {\bf Top}( O^n,  O^n) &
 \stackrel{{\tilde \sigma}_\ast }{\longrightarrow} &{\bf Top}( O^n,
X)\\
 &\tau^{O^n}_n\uparrow & &\uparrow \tau^X_n \\
&
T^n_{O^n}& & T^n_X\\
&\uparrow & & \uparrow \\

&{\bf Costop}(\triangle^n\times{\underline \triangle}, {\underline O^n}) &
&{\bf
Costop}(\triangle^n\times{\underline \triangle}, \underline {X})\\
& \vartheta \uparrow & & \uparrow \vartheta\\
&{\bf Top}(\triangle^n, \vert \vert {\underline  O^n}\vert \vert)&
\stackrel{C_n\vert{\underline {\tilde \sigma}}\vert }\longrightarrow
&{\bf Top}(\triangle^n,\vert  \vert {\underline X}\vert \vert)
\end{array}
$$

\vspace{3mm}
\noindent {\bf 3.10 End of the proof of theorem 2.} It remains to show that
 $
\Theta'_X: BN^*X \rightarrow C^* \vert \vert \underline Y\vert
\vert \cong  C^*\Omega X
$ induces an isomorphism of graded algebras in
(co)homology.  Since it has already been proved that $\Theta'_X$  induces a
linear isomorphism in (co)homology (proposition 3.7), we have only to show that
this isomorphism is compatible with the products.

This is proved  as in lemma 3.7,  by setting:

$F'_1X = BN^*X\otimes BN^\ast X$,  \quad \quad $F'_2X = C^*\vert \vert
\underline Y\vert
\vert
$,
\quad \quad  $\beta'_X = \Theta'_X\circ \nu$ and  $\gamma'_X = \cup
\circ \Theta'_X\otimes \Theta'_X$ where $\nu$ is the product defined in
I-3.3 and $\cup$ is the cup product in
$C^* \vert \vert \underline  Y\vert \vert$ and
considering the models $ {\large \vee}_{p\geq
0}({\large \vee}_0^{p}(\triangle ^n\times\triangle ^p))$ instead of $O^n$.

\vspace{1 cm}

\noindent{\bf  4. Examples.}
\vspace{3mm}

\noindent{\bf 4.1}  $ X = S^{2n}$, $n\geq 1$. For the convenience of the
reader, we give some details
on the computation of $HH^*(X; {\Bbb F}_p)$
when $ X = S^{2n}$, $n\geq 1$ (compare with \cite{[DH2]} and \cite{[KY]}).
First observe that
 as graded vector spaces, $H^\ast \Omega X\cong H_*\Omega X = T(u)$,
$\vert u\vert = 2n-1$, $n\geq 1$. Thus the minimal
 model of $ X = S^{2n}$ is of the form $(TV, d_V)$, as seen in I-$\S
6.1$-a, with
$V\cong s^{-1}T^+(u)$.
More precisely,
$V = \oplus_{k\geq 1} v_k{\Bbb F}_p$ with $v_k$ identified to $s^{-1}u^k$ and
$dv_k = \sum_{i+j=k} v_iv_j$. The map $\mu_V: (T(V'\oplus V''\oplus V'\#V''),
D)\rightarrow (TV, d)$ identifies $V'$ and $V''$ with $V$ and satisfies
$\mu_V(v'_k\# v''_l) ={\footnotesize  \left(\! \! \begin{array}{c} k+l\\l
\end{array} \!\! \right)} v_{k+l}$. For degree reasons, $ X = S^{2n}$ is
${\Bbb F}_p$-$shc$-formal.
On the other hand generators of ${\frak C} _*H^*X$ are of the form
$a_k =  1\underbrace {[ x\vert ...\vert x ]}_{k
\mbox{-times}}$, $k\geq 1$ or
 $b_l =  x\underbrace {[ x\vert ...\vert x ]}_{l
\mbox{-times}}$, $l\geq 0$ when
$H^{2n}X = x{\Bbb F}_p$. The differential is
$ db_l = 0 \,, \quad
da_k  =  \left \{ \begin{array}{ll}
 0, &\mbox{ if } k \mbox{ is odd}\\
2b_{k-1} &\mbox{ if } k \mbox{ is even}
\end{array}
\right.
$. The product in ${\frak C}_*H^*X$ is the composite
 ${\frak C}_*H^*X\otimes {\frak C}_*H^*X
 \stackrel{sh}{\rightarrow} {\frak C}_*(H^*X\otimes H^*X)
\stackrel{{\frak C}_*m}{\rightarrow} {\frak C}_*H^*X$ where $m$ is the
product in $H^*X$. Thus
$$
a_ka_l =   {\footnotesize  \left(\! \! \begin{array}{c} k+l\\l \end{array}
\!\! \right)}
a_{k+l}\,, \quad  b_kb_l = 0\,,\quad
\quad  a_kb_l = {\footnotesize  \left(\! \! \begin{array}{c} k+l\\l
\end{array} \!\! \right)}  b_{k+l}
= b_la_k
$$ if either $k$ or  $l$ is even and the other  products are trivial. We
recover that
$$
H(X^{S^1}; {\Bbb F}_p )=
\left\{
\begin{array}{ll}
\Gamma (a_2)\otimes {\Bbb F}_2[b_0]/b^2_0 \otimes {\Bbb F}_2[a_1]/a_1^2 &
\mbox{ if }\quad  p=2 \\
{\Bbb F}_p\oplus
\oplus_{k\geq 1}a_{2k+1}{\Bbb F}_p\oplus \oplus_{l\geq 0}b_{2l}{\Bbb F}_p
\mbox{ with trivial product}
&\mbox{ if } \quad  p  \mbox{ is odd }
\end{array}
\right.
$$

\vspace{3mm}

\noindent{\bf 4.2}  $ X = {\Bbb C P}(n)$. The space $X$ is $shc$-formal for
every
${\Bbb K}={\Bbb F}_p$, \cite{[Id]}. The graded algebra $H^*(X; {\Bbb F}_p)$
admits a free commutative
minimal model
of the form  $(\Lambda(x, y), d)$ with $dx=0$, \quad
$dy=x^{n+1}$,\quad
 $\vert x\vert =2$ ,\quad  $\vert y\vert=2n+1$. Thus, we deduce from I-$\S
5.3$:
$$ \begin{array}{rl}
  H^*(\Omega X; {\Bbb F}_p) &\cong \Lambda sx\otimes \Gamma sy\\
H^*(X^{S^1}; {\Bbb F}_p) &\cong \left \{
\begin{array}{ll}
 {\Bbb K}[x]/x^{n+1} \otimes \Lambda sx\otimes \Gamma sy
&\mbox{ if }n+1=0 \mbox{ mod } p \\
  \frac{{\Bbb K}[x] \otimes \Lambda sx}{ x^{n+1},(n+1)sx.x^{n+1}} \oplus
\frac{{\Bbb K}^+[x] \otimes \Lambda^+ sx}{ x^{n+1},(n+1)sx.x^{n}}\otimes
\Gamma^+sy
&\mbox{ if not.}
\end{array}
\right.
\end{array}
$$
Observe in particular that the fibration $\Omega X \to X ^{S^1} \to X$ is
not T.N.C.Z.,  if $  n+1
\neq 0$ mod $p$. Compare with
\cite{[MV1]} or \cite{[Ag]}.

\vspace{3mm}

\noindent{\bf 4.3} {\bf  Finite complexes}.  D. Anick proved in \cite{[An]}
that,
if $X$ is  a $r$-connected
CW-complex of dimension
$n$, then,  for any prime $p \geq \frac{n}{r}$, the mod $p$ Adams-Hilton
model is
isomorphic, as a differential Hopf algebra, to the universal enveloping
algebra of
some free differential graded Lie algebra $L $. Therefore, I-5.2 or 5.3 apply.

\vspace{3mm}

\noindent{\bf 4.4}  $ X =\Sigma {\Bbb C P}(2)$ and ${\Bbb K} =
{\Bbb F}_2$. From the  Bott-Samelson theorem, \cite{[Hu]} (Appendix 2),
it is known that $
H_*(\Omega X; {\Bbb F}_2) = T(a_2, b_4)$ with $\triangle b_4 =b_4\otimes 1 +
a_2\otimes a_2 +
1\otimes b_4$. As graded vector spaces,  $H^\ast \Omega X$ is the dual of $
H_\ast \Omega X$
and the product
in $H^\ast \Omega X$ is the dual of the coproduct in $H_\ast \Omega X$.
Thus the element
$a_2$
of degree $2$ in $H^\ast \Omega X$ satisfies $a^2_2=b_4$.
 The minimal model  of $X$ is $(TV, d_V)$ where $V = s^{-1}T^{+}(a_2,
b_4)$. Set
$x_3 = s^{-1}a'_2$, $x_5 = s^{-1}b'_4$.

The map $\mu_V: (T(V'\oplus V''\oplus V'\#V''),
D)\rightarrow (TV, d)$ identifies $V'$ and $V''$ with  $V$.  By 2 of
proposition   I-$\S 6.6$,
 $\mu_0(x'_3\# x''_3) \neq 0$.  (Notice in
passing that the space $ S^3 \vee S^5$ admits the same minimal model as
$X$ but not
the same $shc$-structure). We also have
 $sh(1[ x_3], 1[ x_3]) =  1[
x_3\otimes 1 \vert 1\otimes
x_3]  +
  1[ 1\otimes x_3\vert x_3\otimes 1]$.

 If $s_V$ is any section of ${\frak C}_* (\Psi)$, a
straightforward computation shows
that   $s_V( 1[ x_3\otimes 1\vert
x_3]) =
1[ x'_3\vert  x''_3]$ and
$s_V(1\otimes 1[ 1\otimes x_3\vert x_3\otimes 1] = 1[
x''_3\vert
x'_3] + 1[ x'_3\# x''_3] $, so  that ${\frak C}_\ast \mu_V(sh(1[ x_3], 1[
x_3])) = 1[
x_5]$. Thus, the cohomology class
$[1[ x_3]]^2 = [1[ x_5] ]$  is not trivial
in  $H^*\left( \left( \Sigma {\Bbb C}P\left( p \right)\right)^{S^1}; {\Bbb
F}_p\right)$. More
generally, for
any prime $p$, if $[1[ x_3] ]$ is
a generator of $H^2\left( \left( \Sigma {\Bbb C}P\left( p
\right)\right)^{S^1}; {\Bbb F}_p\right)$,  then
$[1[ x_3]]^p \neq 0.$

\vspace{3mm}

\noindent{\bf 4.5} $X= G_2$ and ${\Bbb K}={\Bbb F}_5$.  From, \cite{[MT]}
( 5.9  and  5.18),  we know that if $X$ is   the exceptional Lie group
$G_2$ then
 $H^*(G_2; {\Bbb F}_5) = \Lambda (v_3, v_{11})$
 is the free commutative algebra
on two generators of degree 3 and 11 respectively, and
the first Steenrod operation satisfies: ${\cal P}^1 v_3 \not = 0$. The Kudo
transgression theorem
implies that there exists $ y_2 \in H^*(\Omega G_2;
{\Bbb F}_5)$ such that $y_2^5 \not = 0$. Since
 the multiplicative fibration
$\Omega G_2 \to G_2^{S^1} \to G_2$ admits a section, $G_2^{S^1} \simeq G_2
\times \Omega G_2$. Thus there is a cohomology class $\alpha \in H^{2} (
G_2^{S^1}; {\Bbb F}_5)$ such that $\alpha ^5 \not = 0$. We recover this result
with our construction.

The minimal model
of $G_2$ is given by $(TV, d_V)$ where $V \cong s^{-1} H^*(\Omega G_2; {\Bbb
F}_5) $.
 More
precisely, this minimal model is $(T(x_3, x_5, x_7, x_9, x_{11}, y_{11},...),
d)$ with $dx_3 = dy_{11} = 0$, $dx_5 = x^2_3$,
$dx_7 = [x_3, x_5]$, $dx_9 = [x_3, x_7] + x^2_5$ and $dx_{11} = [x_3, x_9]
>+ [x_5, x_7]$.\\
Consider the map   $ \mu_V: (T(V'\oplus V''\oplus V'\#V''),
D)\rightarrow (TV, d)$. We have $\mu_V(x'_3\# x''_3) =
2x_5$; $\mu_V(x'_3\# x''_5) = 3x_7$;
$\mu_V(x'_3\# x''_7) = 4x_9$; $\mu_V(x'_3\# x''_9) =
\epsilon y_{11}$ with $\epsilon \neq 0$,  by 2 of proposition I-$\S 6.6$. A
section
 $s_V$ of ${\frak C}_* (\Psi)$ can be determined in low degrees and  a
tedious but straightforward
computation shows that the cocycle
$(1 [ x_3 ]  )^5 $ is cohomologous to $4\epsilon (1[ y_{11} ] )$ and thus
non trivial. By contrast,
the cohomology  class of $(1[ x_3 ])^{25}$ is zero as  can be shown
by a straightforward
  computation.

\vspace{0,5 cm}

       \end{document}